\documentclass{preprint}

\usepackage{amsmath,amsfonts,amssymb,amsthm,mathtools}
\usepackage{hhline}

\makeatletter
\DeclareRobustCommand*\cal{\@fontswitch\relax\mathcal}
\makeatother

\newcommand{\E}[1]{\mathop{{\rm \bf E}\!\left\{#1\right\}}\nolimits}

\newtheorem{rem}{Remark}
\newtheorem{exmp}{Example}

\begin{document}

\title{Improved Discrete-Time Kalman Filtering within Singular Value Decomposition}

\author{\au{Maria V. Kulikova$^{1}$}, 
\au{Julia V. Tsyganova$^2$}
}

\address{\add{1}{CEMAT (Center for Computational and Stochastic Mathematics), Instituto Superior T\'ecnico, Universidade de Lisboa, Av.~Rovisco Pais, 1049-001 Lisboa, Portugal}
\add{2}{Department of Mathematics, Information and Aviation Technology, Ulyanovsk State University, Str. L. Tolstoy 42, 432017 Ulyanovsk, Russia}
\email{maria.kulikova@ist.utl.pt; tsyganovajv@gmail.com}}

\begin{abstract}
The paper presents a new Kalman filter (KF) implementation useful in applications where the accuracy of numerical solution of the associated Riccati equation might be crucially reduced by  influence of roundoff  errors. Since the appearance of the KF in 1960s, it has been recognized that the factored-form of the KF is preferable for practical implementation. The most popular and beneficial techniques are found in the class of square-root algorithms based on the Cholesky decomposition of error covariance matrix. Another important matrix factorization method is the singular value decomposition (SVD) and, hence, further encouraging implementations might be found under this approach. The analysis presented here exposes that the previously proposed SVD-based KF variant is still sensitive to roundoff errors and poorly treats ill-conditioned situations, although the SVD-based strategy is inherently more stable than the conventional KF approach. In this paper we design a new SVD-based KF implementation for enhancing the robustness against roundoff errors, provide its detailed derivation, and discuss the numerical stability issues. A set of numerical experiments are performed for comparative study. The obtained results illustrate that the new SVD-based method is algebraically equivalent to the conventional KF and to the previously proposed SVD-based method, but it outperforms the mentioned techniques for estimation accuracy in ill-conditioned situations.
\end{abstract}

\maketitle


\section{Introduction}

The main issue to be addressed in this paper is a numerical robustness of the Kalman filter (KF). It is known that in some applications, the accumulation of roundoff errors might severely affect computations and eventually degrade the performance of the filter~\cite{VerhaegenDooren1986,Verhaegen1989}. Due to this fact, since 1960s special attention has been paid in the KF community for designing robust filter implementations~\cite[p.~2058]{Grewal2010}: ``It was recognized during the early years of Kalman filter applications that factored-form Kalman filters (square root filters) are the preferred implementation for applications demanding high operational reliability''. The key idea of the factored-form KF strategy is to factorize the error covariance matrix $P$ in the form of $P = SS^T$ and, then, re-formulate the KF equations in terms of these factors only. Thus, only the resulted factors are processed between time and measurement updates in each iteration step of the KF recursion. Undoubtedly, this approach is not free of roundoff errors, however, the square-root (SR) methodology ensures that the covariance matrix, reformed by multiplying its factors (i.e. $SS^T=P$), is symmetric and positive semi-definite~\cite[Chapter~7]{GrewalAndrews2015}.
Additionally, the SR approach allows for calculation with double precision~\cite[pp.~728-729]{KaminskiBryson1971}: ``The advantage of the SR approach becomes obvious by relating the condition number $K(P)$ and $K(S)$, i.e. $K(P) = K(SS^T) = [K(S)]^2$. This means that while numerical operation with $P$ may encounter difficulties when $K(P) = 10^p$, the SR filter should function until $K(P) = 10^{2p}$, i.e. with double precision.''
Additionally, the modern KF algorithms are formulated in the so-called {\it array} form that implies utilization of numerically stable orthogonal transformations for updating the covariance matrix factors~\cite[Chapter~12]{KailathSayed2000}.
More precisely, orthogonal operators are applied to the {\it pre-array} (which accommodates the filter quantities available at the current step) in order to get a special form of the {\it post-array}. Next, the updated filter quantities are simply read off from the post-array. This feature makes the array algorithms better suited to parallel implementation and to very large scale integration (VLSI) implementation as mentioned in~\cite{ParkKailath1995}. Besides, they are easier to implement than the explicit filter equations. From a numerical stability standpoint, the use of orthogonal transformations improves a reliability of the estimator; see, for instance, the discussion in~\cite{KuKu2017IET_CTA} and many others.

We stress that the covariance square roots are not uniquely determined by relationship $P=SS^T$. As discussed in~\cite[p.~958]{Haykin2007}: ``There are many types of matrix decomposition techniques that factorize a covariance matrix P in the form of $P = SS^T$. E.g., the Cholesky decomposition, the singular value decomposition (SVD) and the eigenvector decomposition. Choosing which one to use primarily depends upon the particular application, numerical concerns and desired level of accuracy.''

The SR approach based on the Cholesky decomposition is the most popular technique for enhancing the numerical robustness of the filter; see~\cite{KaminskiBryson1971,Morf1975,Park1994a,Sayed1994} and many others. It implies the error covariance matrix decomposition in the form of $P = P^{1/2}P^{T/2}$ where $P^{1/2}$ is an upper or a lower triangular matrix with positive diagonal elements, i.e. $S=P^{1/2}$. The important fact to be mentioned here is that the Cholesky decomposition  exists and is unique when the symmetric matrix to be decomposed is positive definite~\cite{Golub1983}. If the matrix is a positive semi-definite, then the Cholesky decomposition still exists, however, it is not unique~\cite{Higham1990}.

The UD factorization-based strategy developed in~\cite{Bierman1977} is a `square-root-free' modification of the factored-form KF approach based on the Cholesky decomposition. The UD-based methods imply the decomposition of the error covariance matrix in the form of $P = \bar U D \bar U^T$ where $\bar U$ is an upper triangular or a lower triangular matrix with 1's on the main diagonal and $D$ is a diagonal matrix. Hence, in case of UD-based KF algorithms, the square-root factor of $P$ can be chosen to be $S=\bar U D^{1/2}$. In the UD-based class of methods, the modified weighted Gram-Schmidt (MWGS) orthogonalization (see~\cite{Bjorck1967}) is used for recursive update of the resulted factors $\bar U$ and $D$ in each iteration step of the filter. The UD-based KF implementations were shown to have roughly the same computational complexity as the conventional KF implementation, but with better numerical stability as in all factorized-form (SR) filter implementations~\cite{bierman1977numerical}.

Alternatively, in this paper we explore the SVD factorization-based KF approach. It implies the decomposition of the error covariance matrix in the form of $P = Q \Sigma Q^T$ where $Q$ is an orthogonal matrix and $\Sigma$ is a diagonal matrix with singular values of $P$. Hence, in case of SVD-based KF algorithms, the square-root factor of $P$ can be chosen to be $S=Q \Sigma^{1/2}$. The SVD is one of the most important and useful matrix decompositions in linear algebra and it is the most robust algorithm to factorize a covariance matrix especially when it becomes nearly singular~\cite[p.~2596]{Menegaz2015}: ``From a numerical implementation standpoint, even though the Cholesky decomposition seems to be the most adopted method to compute the square-root matrix of the state covariance, some studies indicate that other methods, such as SVD decomposition, provide better estimation quality''. Some evidences of better estimation quality exist in the field of nonlinear filtering~\cite{zhang2008,straka2013} {\it etc}. For linear filtering problem
examined in this paper, the first SVD-based KF was, to the best of
our knowledge, developed in~\cite{WangSVD1992}. In this paper, we explain how it can be improved further for enhancing its numerical robustness against roundoff errors. More precisely, we design a new SVD-based KF implementation, provide its detailed derivation and discuss the numerical stability issues. The filter is derived here in terms of the {\it covariance} quantities and it is formulated in the {\it array} form. A set of numerical experiments are performed for comparative study.

Finally, we would like to note that the filter developed in this paper has been recently used in adaptive gradient-based KF schemes~\cite{TsyganovaOnline}. In the cited work, the UD- and the SVD-based techniques were shown to provide the best estimation quality when solving ill-conditioned parameter estimation problems. This creates a strong background for their practical use. Thus, the potential applications of the SVD-based methodology developed in this paper include, e.g. the space navigation problem for orbit determination case study~\cite{bierman1977numerical,Hall1996}, the calibration and alignment of inertial navigation systems~\cite{Grewal1991}, the development of the Global Positioning Systems (GPS)~\cite{Grewal2007} and many other applications discussed in~\cite{Grewal2010a}.

\section{Conventional Kalman filter implementation}

Consider the state-space equations
  \begin{align}
   x_{k}& =  F_{k-1} x_{k-1}+B_{k-1}u_{k-1}+ G_{k-1}w_{k-1}, && k \ge 1,  \label{eq:st:1} \\
    z_k & =  H_k x_k+v_k &&              \label{eq:st:2}
  \end{align}
where $x_k \in \mathbb R^n$, $u_k \in \mathbb R^d$ and $z_k \in \mathbb R^m$ are,
respectively, the unknown dynamic state, the deterministic control input and the available measurement vectors. The process and the measurement noises are independent Gaussian zero-mean white-noise processes with covariances $\Theta_k \ge 0$ and $R_k > 0$, respectively, i.e. $w_k \sim {\mathcal N}\left(0,\Theta_k\right)$ and $v_k  \sim {\mathcal N}\left(0,R_k\right)$.  They are also uncorrelated with the initial state $x_0 \sim {\mathcal N}\left(\bar x_0,\Pi_0\right)$, $\Pi_0 \ge 0$.

The KF yields the linear minimum least-square estimate of the dynamic state, $\hat x_{k|k}$, given the measurements $\{z_1,\ldots, z_{k}\}$, which can be recursively computed as follows~\cite[Theorem~9.2.1]{KailathSayed2000}:

\textsc{Initialization}  ($k=0$). Set $\hat x_{0|0} = \bar x_0$ and $P_{0|0} = \Pi_0$.

\textsc{Time Update} ($k=1, \ldots, K$). The step of going from $\hat x_{k-1|k-1}$ to $\hat x_{k|k-1}$ is called the time update. At this stage, the one-step ahead predicted ({\it a priori}) estimate, $\hat x_{k|k-1}$, is computed together with the corresponding error covariance matrix $P_{k|k-1}=\E{ (x_{k}-\hat x_{k|k-1})(x_{k}-\hat x_{k|k-1})^T}$ as follows:
\begin{align}
\hat x_{k|k-1} & = F_{k-1}\hat x_{k-1|k-1}+B_{k-1}u_{k-1}, \label{kf:p:X} \\
P_{k|k-1} & = F_{k-1}P_{k-1|k-1}F_{k-1}^T+G_{k-1}\Theta_{k-1}G_{k-1}^T. \label{kf:p:P}
\end{align}

\textsc{Measurement Update}  ($k=1, \ldots, K$). When the new measurement $z_k$ becomes available,  the recent information is extracted for updating the state estimate through the feedback gain $K_k$. This correction step is called the measurement update where {\it a posteriori} estimate $\hat x_{k|k}$ is calculated together with the corresponding error covariance matrix $P_{k|k}$ as follows:
\begin{align}
K_{k} & = P_{k|k-1}H^T_kR_{e,k}^{-1}, \!\!\!  & R_{e,k} & = H_kP_{k|k-1}H_k^T+R_k, \label{kf:f:K} \\
\hat x_{k|k} & =    \hat x_{k|k-1}+K_{k}e_k,\! \! \!  & e_k & = z_k-H_k\hat x_{k|k-1},   \label{kf:f:X} \\
P_{k|k}  & = (I - K_{k}H_k)P_{k|k-1}. \! \! \! \! && \label{kf:f:P}
\end{align}
The important property of the KF for Gaussian state-space models~\eqref{eq:st:1}, \eqref{eq:st:2}
is $e_k \sim {\cal N}\left(0, R_{e,k}\right)$ where $\{ e_k \}$ are innovations.

For further discussion of the KF implementation methods, some auxiliary theoretical results are required.  More precisely, the following formulas hold for the classical KF~\cite[p.~128-129]{simon2006optimal}:
\begin{align}
\!\!K_{k} & = P_{k|k-1}H_k^T(H_kP_{k|k-1}H_k^T+R_k)^{-1} \label{kf:K:eq1} \\
 & = P_{k|k}H_k^TR_k^{-1}. \label{kf:K:eq2} \\
\!\!P_{k|k}  & = \left( P_{k|k-1}^{-1}+H^T_kR_k^{-1}H_k\right)^{-1} \label{kf:P:eq1} \\
  & = (I - K_{k}H_k)P_{k|k-1} \label{kf:P:eq2} \\
  & = (I - K_{k}H_k)P_{k|k-1}(I - K_{k}H_k)^T+K_kR_kK_k^T. \label{kf:P:eq3}
\end{align}

\section{The SVD factorization-based KF} \label{sec:svd}

Consider the SVD factorization~\cite[Theorem~1.1.6]{Bjorck2015}: Every matrix $A \in {\mathbb C}^{m\times n}$ of rank $r$ can be written as
\[
A  = W\Sigma V^*, \,
\Sigma =
\begin{bmatrix}
S & 0 \\
0 & 0
\end{bmatrix} \in {\mathbb C}^{m\times n},  \; S={\rm diag}\{ \sigma_1,\ldots,\sigma_r\}
\]
where $W \in {\mathbb C}^{m\times m}$, $V \in {\mathbb C}^{n\times n}$  are unitary matrices, $V^*$ is the conjugate transpose of $V$, and $S \in {\mathbb R}^{r\times r}$ is a real nonnegative diagonal matrix. Here $\sigma_1\geq \sigma_2\geq \ldots \geq\sigma_r>0$ are called the
singular values of $A$. (Note that if $r = n$ and/or $r = m$, some of the zero submatrices in $\Sigma$ are empty.)



In the KF realm, the preferred factorization method is applied to the state error covariance matrix involved in equations~\eqref{kf:p:X}~-- \eqref{kf:f:P}. More precisely, the SVD factorization is performed only once, i.e. for the initial error covariance $\Pi_0 \in {\mathbb R}^{n\times n}$. We have $\Pi_0 = Q_{\Pi_0}D_{\Pi_0}Q_{\Pi_0}^T$  where $Q_{\Pi_0}$  and $D_{\Pi_0}$ are an orthogonal and a diagonal matrices, respectively. The matrix $D_{\Pi_0}$ contains the singular values of $\Pi_0$. Next, the conventional KF equations~\eqref{kf:p:X}~-- \eqref{kf:f:P}  are re-formulated in terms of the resulted factors only. Thus, the SVD-based KF algorithms recursively update the SVD factors $D_{P_{k|k}}$ and $Q_{P_{k|k}}$ of the error covariance matrix $P_{k|k}$ in each iteration step of the filter instead of calculating the full matrix $P_{k|k}$. The previous works on the SVD factorization-based KF have produced the following implementation; see Eqs~(17), (22), (23) in~\cite[pp.~1226]{WangSVD1992}:

\textsc{Initialization}  ($k=0$). Apply SVD factorization for the initial
error covariance matrix $\Pi_0 = Q_{\Pi_0}D_{\Pi_0}Q_{\Pi_0}^T$. Set the initial values as follows: $Q_{P_{0|0}} = Q_{\Pi_0}$, $D^{1/2}_{P_{0|0}} = D^{1/2}_{\Pi_0}$ and $\hat x_{0|0} = \bar x_0$.

\textsc{Time Update} ($k=1, \ldots, K$). Given $\hat x_{k-1|k-1}$, compute {\it  a priori} estimate $\hat x_{k|k-1}$ by equation~\eqref{kf:p:X} and find the SVD factors, $\{Q_{P_{k|k-1}}, D^{1/2}_{P_{k|k-1}}\}$, of the corresponding $P_{k|k-1}$ as follows:
\begin{align}
\underbrace{
\begin{bmatrix}
D^{1/2}_{P_{k-1|k-1}}Q^T_{P_{k-1|k-1}}F_{k-1}^T\\
\Theta_{k-1}^{1/2}G_{k-1}^T
\end{bmatrix}
}_{\rm Pre-array}
 & =
 \underbrace{
\mathfrak{W}_{TU}
\begin{bmatrix}
D_{P_{k|k-1}}^{1/2} \\
0
\end{bmatrix}
\mathfrak{V}^T_{TU},
}_{\rm Post-array \: SVD \:factors}    \label{svd:p:fac} \\
Q_{P_{k|k-1}} & = \mathfrak{V}_{TU} \label{svd:p:P}
\end{align}
where $\Theta_{k-1}^{1/2}$ is a Cholesky factor of the process noise covariance $\Theta_{k-1}$ such that $\Theta_{k-1} = \Theta_{k-1}^{T/2}\Theta_{k-1}^{1/2}$, where $\Theta_{k-1}^{1/2}$ is an upper triangular matrix with positive diagonal elements. The orthogonal matrices $\mathfrak{W}_{TU}$, $\mathfrak{V}_{TU}$ are the resulted post-array SVD factors.

\textsc{Measurement Update}  ($k=1, \ldots, K$). Given $\hat x_{k|k-1}$, compute {\it a posteriori} estimate $\hat x_{k|k}$ by equation~\eqref{kf:f:X}. Next, given the SVD factors $\{Q_{P_{k|k-1}}, D^{1/2}_{P_{k|k-1}}\}$, calculate their {\it a posteriori} values $\{Q_{P_{k|k}}, D^{1/2}_{P_{k|k}}\}$ as follows:
\begin{align}
\underbrace{
\begin{bmatrix}
R_k^{-T/2}H_kQ_{P_{k|k-1}} \\
D_{P_{k|k-1}}^{-1/2}
\end{bmatrix}
}_{\mbox{\small Pre-array}}
 & =
\underbrace{
\mathfrak{W}_{MU}
\begin{bmatrix}
D_{P_{k|k}}^{-1/2} \\
0
\end{bmatrix}
\mathfrak{V}_{MU}^T,
}_{\mbox{\small Post-array SVD factors}} \label{svd:f:fac} \\
Q_{P_{k|k}} & =  Q_{P_{k|k-1}}\mathfrak{V}_{MU}, \label{svd:f:P} \\
K_{k} & = \left(Q_{P_{k|k}}D_{P_{k|k}}Q^T_{P_{k|k}}\right)H_k^TR_k^{-1} \label{svd:f:K}
\end{align}
where $R_k^{1/2}$ is a Cholesky factor of the measurement covariance $R_k$ such that $R_k = R_k^{T/2}R_k^{1/2}$, where $R_k^{1/2}$  is an upper triangular matrix with positive diagonal elements. Additionally, $R_k^{-T/2} = (R_k^{T/2})^{-1}$. The orthogonal matrices $\mathfrak{W}_{MU}$, $\mathfrak{V}_{MU}$ are the post-array SVD factors.

\begin{rem} \label{rem:1}
Each iteration step of the examined SVD-based filter  has the following pattern: given a pre-array $A \in {\mathbb R}^{(k+s)\times s}$, compute the post-array SVD factors $\mathfrak{W} \in {\mathbb R}^{(k+s)\times (k+s)}$, $\Sigma \in {\mathbb R}^{(k+s)\times s}$ and $\mathfrak{V} \in {\mathbb R}^{s\times s}$ by means of the SVD factorization, $A=\mathfrak{W}\Sigma\mathfrak{V}^T$.
\end{rem}

The readers are referred to~\cite{WangSVD1992} for the proof and detailed explanation of algebraical equivalence between the conventional KF~\eqref{kf:p:X}~-- \eqref{kf:f:P} and the SVD factorization-based version presented above. Here, we stress that the SVD-based KF variant implies formula~\eqref{kf:P:eq1} for computing {\it a posteriori} error covariance matrix and equation~\eqref{kf:K:eq2} for calculating the feedback gain $K_k$. It is easy to prove that if we take into account the SVD factorization $A=\mathfrak{W} \Sigma \mathfrak{V}^T$ and the properties of orthogonal matrices. Indeed, for each pre-array to be decomposed we have $A^TA = (\mathfrak{V} \Sigma \mathfrak{W}^T)(\mathfrak{W}\Sigma \mathfrak{V}^T) = \mathfrak{V} \Sigma^2 \mathfrak{V}^T$. Next, by comparing both sides of the obtained matrix equations, we come to the corresponding SVD-based KF formulas.

To begin constructing the new SVD-based implementation, we first discuss the filter presented above. Our analysis exposes that two $n \times n$  matrices' inversions are required in each iteration step of the examined implementation; see the terms $D_{P_{k|k-1}}^{-1/2}$ and $D_{P_{k|k}}^{-1/2}$ in equation~\eqref{svd:f:fac} and, next, compare this equation with formula~\eqref{svd:p:fac} where only $D_{P_{k|k-1}}^{1/2}$ is available from the post-array and $D_{P_{k|k}}^{1/2}$ is required in the pre-array for performing the next step of the recursion. The matrices to be inverted are diagonal and, hence, the inversion is simply $n$ scalar divisions by the corresponding diagonal entries. However, in the SVD-based filtering methods, these diagonal elements are, in fact, the square roots of non-zero sigma values of the state error covariances. If the error covariance matrix is close to singular, then the mentioned matrix inversion implies the scalar divisions by very small numbers that causes large roundoff errors or even obstructs the inversion. Hence, we anticipate that the influence of roundoff errors deteriorates the estimation accuracy of the previously developed SVD-based KF variant rapidly while solving ill-conditioned problems. The results of numerical experiments presented in Section~\ref{sec:exp} substantiate our theoretical expectations.

Finally, the filter implementation summarized by equations~\eqref{kf:p:X}, \eqref{kf:f:X}, \eqref{svd:p:fac}~-- \eqref{svd:f:K}  requires the Cholesky decomposition of the noise covariance matrices $\Theta_k$ and $R_k$. Hence, we further abbreviate this algorithm as \texttt{SVD-SRKF}, emphasizing the requirement of the Cholesky factors in this implementation. Evidently, for time-invariant matrices $\Theta$ and $R$ the Cholesky decomposition is performed only once. This can be done at the initial step, i.e.  before the KF recursion. It should be also pointed out that the \texttt{SVD-SRKF} additionally requires the $m \times m$ matrix  inversion; see  the terms $R^{-T/2}_k$ and $R_k^{-1}$ in equations~\eqref{svd:f:fac}, \eqref{svd:f:K}. Again, if the process and measurement noise covariances are constant over time, then the inverse matrices might be pre-computed while initializing the filter.

\section{Main result} \label{sec:main}

To enhance the estimation accuracy of the original SVD-based KF and, hence, to increase the practical applicability of the SVD factorization-based KF strategy, we derive a new implementation under this approach. Recall that the original \texttt{SVD-SRKF} requires two $n \times n$ diagonal matrices' inversions as explained in Section~\ref{sec:svd}. Our first suggestion for the filter improvement is to utilize equation~\eqref{kf:P:eq3} for computing {\it a posteriori} error covariance $P_{k|k}$ instead of~\eqref{kf:P:eq1}. Such symmetric form is known as the Joseph stabilized implementation and allows for reducing the influence of roundoff errors as discussed in~\cite{VerhaegenDooren1986}, \cite{Verhaegen1989}. It also ensures the error covariance matrix $P_{k|k}$ to be symmetric. Besides, it guarantees the positive definiteness of $P_{k|k}$ as long as $P_{k|k-1}$ is positive definite matrix~\cite[p.~129]{simon2006optimal}. Hence, our goal is to derive the SVD-based KF algorithm based on equation~\eqref{kf:P:eq3}.


Our next suggestion concerns the Cholesky decomposition of the process and measurement noise covariances $\Theta_k$ and $R_k$ required in the previously-proposed \texttt{SVD-SRKF} implementation. We stress that in order to perform the Cholesky decomposition, the mentioned symmetric matrices should be positive-definite. If any of them becomes positive semi-definite in any recursion step, then this operation causes an unexpected failure of the corresponding practical implementation, because for positive semi-definite matrix the Cholesky decomposition exists but not unique. When implementing in MATLAB, the code will be interrupted by the error appeared in function \verb"chol", which performs the Cholesky decomposition. Hence, we formulate the method that utilizes the SVD matrix factorization only. The new implementation is abbreviated further to \texttt{SVD-KF}, emphasizing the requirement of only SVD factors.

Now, we are ready to present the new SVD-based KF implementation. Instead of conventional recursion~\eqref{kf:p:X}~-- \eqref{kf:f:P} for $P_{k|k}$, we update its SVD factors, $\{Q_{P_{k|k}}, D^{1/2}_{P_{k|k}}\}$, as follows:

\textsc{Initial Step}  ($k=0$). Apply SVD factorization for the initial error covariance matrix $\Pi_0 = Q_{\Pi_0}D_{\Pi_0}Q_{\Pi_0}^T$. Set the initial values as follows: $Q_{P_{0|0}} = Q_{\Pi_0}$, $D^{1/2}_{P_{0|0}} = D^{1/2}_{\Pi_0}$ and $\hat x_{0|0} = \bar x_0$.

\textsc{Time Update} ($k=1, \ldots, K$).  Given $\hat x_{k-1|k-1}$, compute {\it  a priori} estimate $\hat x_{k|k-1}$ by equation~\eqref{kf:p:X}. Build the pre-array and apply the SVD factorization in order to obtain {\it a priori} error covariance SVD factors $\{Q_{P_{k|k-1}}, D_{P_{k|k-1}}^{1/2}\}$ as follows:
\begin{align}
\!\!\!
\underbrace{
\begin{bmatrix}
D^{1/2}_{P_{k-1|k-1}}Q^T_{P_{k-1|k-1}}F_{k-1}^T\\
D^{1/2}_{\Theta_{k-1}} Q^T_{\Theta_{k-1}} G_{k-1}^T
\end{bmatrix}
}_{\rm Pre-array} \!\!
=
\underbrace{
\mathfrak{W}_{TU}
\begin{bmatrix}
D_{P_{k|k-1}}^{1/2} \\
0
\end{bmatrix}
Q_{P_{k|k-1}}^T
}_{\rm Post-array \: SVD \:factors} \!\!\! \label{svd:5}
\end{align}
where  $\{ Q_{\Theta_{k-1}}, D_{\Theta_{k-1}} \}$ are the SVD factors of the process noise covariance $\Theta_{k-1}$, i.e. $\Theta_{k-1} = Q_{\Theta_{k-1}}D_{\Theta_{k-1}}Q_{\Theta_{k-1}}^T$.

\textsc{Measurement Update}  ($k=1, \ldots, K$). Build the pre-arrays from the filter quantities that are currently available and, then, apply the SVD factorizations in order to obtain the corresponding SVD factors of the updated filter quantities as follows:
\begin{align}
\underbrace{
\begin{bmatrix}
D^{1/2}_{R_k} Q^T_{R_k} \\
D^{1/2}_{P_{k|k-1}}Q^T_{P_{k|k-1}}H_k^T
\end{bmatrix}
}_{\rm Pre-array}
 =
\underbrace{
\mathfrak{W}_{MU}^{(1)}
\begin{bmatrix}
D_{R_{e,k}}^{1/2} \label{svd:1} \\
0
\end{bmatrix}
Q_{R_{e,k}}^T
}_{\rm Post-array \: SVD \:factors},  \\
\bar K_{k} = P_{k|k-1}H_k^TQ_{R_{e,k}}, \quad K_k = \bar K_k D^{-1}_{R_{e,k}}Q^T_{R_{e,k}}, \label{svd:K} \\
\!\!\!\!\!\!\!
\underbrace{
\begin{bmatrix}
D_{P_{k|k-1}}^{1/2}Q_{P_{k|k-1}}^T\left(I - K_k H_k\right)^T \\
D^{1/2}_{R_k} Q^T_{R_k}  K_{k}^T
\end{bmatrix}
}_{\rm Pre-array}
\!\! =
\underbrace{
\mathfrak{W}_{MU}^{(2)}
\begin{bmatrix}
D_{P_{k|k}}^{1/2} \\
0
\end{bmatrix}
Q_{P_{k|k}}^T
}_{\rm Post-array \: SVD \:factors} \label{svd:2}
\end{align}
where $\{ Q_{R_{k}}, D_{R_{k}} \}$ are the SVD factors of the measurement covariance $R_k$, i.e. $R_k = Q_{R_k}D_{R_k}Q_{R_k}^T$. The matrices $\mathfrak{W}_{MU}^{(1)} \in {\mathbb R}^{(m+n)\times (m+n)}$, $Q_{R_{e,k}} \in {\mathbb R}^{m\times m}$ and $\mathfrak{W}_{MU}^{(2)} \in {\mathbb R}^{(n+m)\times (n+m)}$, $Q_{P_{k|k}} \in {\mathbb R}^{n\times n}$ are the orthogonal matrices of the corresponding SVD factorizations in~\eqref{svd:1} and \eqref{svd:2}, respectively. Next, $D_{R_{e,k}}^{1/2} \in {\mathbb R}^{m\times m}$ and $D_{P_{k|k}}^{1/2} \in {\mathbb R}^{n\times n}$ are the diagonal matrices with square roots of the singular values of $R_{e,k}$ and $P_{k|k}$, respectively.
Finally, find {\it a posteriori} estimate $\hat x_{k|k}$ through equations
\begin{align}
\bar e_k     & =  Q_{R_{e,k}}^T e_k, \quad e_k = z_k-H_k\hat x_{k|k-1}, \\
\hat x_{k|k} & = \hat x_{k|k-1} + \bar K_{k}D^{-1}_{R_{e,k}}\bar e_k,   \label{svd:3}
\end{align}
where the required SVD factors of the innovation covariance $R_{e,k}$, i.e. $\{Q_{R_{e,k}}, D_{R_{e,k}}\}$, are directly read-off from the post-array factors in equation~\eqref{svd:1}.

To validate the new method presented above, we prove the algebraic equivalence between the \texttt{SVD-KF} recursion~\eqref{kf:p:X},~\eqref{svd:5}~-- \eqref{svd:3}  and the conventional KF~\eqref{kf:p:X}~-- \eqref{kf:f:P}. From a general form of the SVD-based filter iterates (i.e. $A=\mathfrak{W} \Sigma \mathfrak{V}^T$), we obtain  $A^TA = (\mathfrak{V} \Sigma \mathfrak{W}^T)(\mathfrak{W}\Sigma \mathfrak{V}^T) = \mathfrak{V} \Sigma^2 \mathfrak{V}^T$. Next, comparing both sides of the resulted matrix equality $A^TA=\mathfrak{V} \Sigma^2 \mathfrak{V}^T$ we derive the required formulas. More precisely, from~\eqref{svd:5} of the \texttt{SVD-KF} we obtain
\begin{align*}
P_{k|k-1} & = Q_{P_{k|k-1}}D_{P_{k|k-1}}Q_{P_{k|k-1}}^T  \\
& = F_{k-1}Q_{P_{k-1|k-1}}D_{P_{k-1|k-1}}Q_{P_{k-1|k-1}}^T F_{k-1}^T \\
&+ G_{k-1}Q_{\Theta_{k-1}}D_{\Theta_{k-1}} Q^T_{\Theta_{k-1}}G_{k-1}^T \\
& = F_{k-1}P_{k-1|k-1}F_{k-1}^T+G_{k-1}\Theta_{k-1}G_{k-1}^T,
\end{align*}
which is exactly formula~\eqref{kf:p:P}  of the conventional KF.

Next, we validate formula~\eqref{svd:1} of the \texttt{SVD-KF}. Taking into account the properties of orthogonal matrices, from~\eqref{svd:1} we obtain
\begin{align*}
R_{e,k} & = Q_{R_{e,k}}D_{R_{e,k}}Q_{R_{e,k}}^T  \\
        & = Q_{R_k}D_{R_k}Q_{R_k}^T + H_k Q_{P_{k|k-1}} D_{P_{k|k-1}}Q^T_{P_{k|k-1}}H_k^T \\
        & = R_k + H_kP_{k|k-1}H_k^T,
\end{align*}
i.e. formula~\eqref{svd:1} implies $R_{e,k} = H_kP_{k|k-1}H_k^T+R_k$ that is exactly the second formula in equation~\eqref{kf:f:K} of the conventional KF.

Next, expression~\eqref{svd:K} for calculating the feedback gain $K_k$  and its normalized variant $\bar K_k$ is derived from~\eqref{kf:K:eq1} where the matrix $R_{e,k}$ is the SVD factorized. Indeed,
\begin{align*}
K_k  & = P_{k|k-1}H_k^T(H_kP_{k|k-1}H_k^T+R_k)^{-1} = P_{k|k-1}H_k^TR_{e,k}^{-1} \\
& = P_{k|k-1}H^T_k \left(Q_{R_{e,k}} D_{R_{e,k}}Q_{R_{e,k}}^T\right)^{-1} \\
& = P_{k|k-1}H^T_k Q_{R_{e,k}}D_{R_{e,k}}^{-1}Q_{R_{e,k}}^T = \bar K_k D_{R_{e,k}}^{-1}Q_{R_{e,k}}^T
\end{align*}
where we have introduced an expression for the normalized feedback gain as follows: $\bar K_k  = P_{k|k-1}H^T_k Q_{R_{e,k}}$.

The algebraic equivalence between~\eqref{kf:P:eq3} and~\eqref{svd:2} can be proved at the same manner as it was shown for $R_{e,k}$, i.e. from~\eqref{svd:2} we obtain
\begin{align*}
P_{k|k} & = Q_{P_{k|k}}D_{P_{k|k}}Q^T_{P_{k|k}}\\
& = \left(I - K_{k} H_k\right)Q_{P_{k|k-1}}D_{P_{k|k-1}}Q^T_{P_{k|k-1}} \left(I - K_{k}  H_k\right)^T \\
& \phantom{=}\; + K_k Q_{R_k} D_{R_k}Q_{R_k}^{T} K_k^T.
\end{align*}

Finally, equation~\eqref{svd:3} for computing {\it a posteriori} state estimate is derived from~\eqref{kf:f:X} as follows:
\begin{align*}
\hat x_{k|k} & = \hat x_{k|k-1}+K_{k}e_k = \hat x_{k|k-1}+\bar K_k D_{R_{e,k}}^{-1}Q_{R_{e,k}}^T e_k \\
&  = \hat x_{k|k-1}+\bar K_{k} D_{R_{e,k}}^{-1} \bar e_k \quad {\rm where} \quad \bar e_k = Q_{R_{e,k}}^T e_k.
\end{align*}
This concludes the proof. Thus, the new \texttt{SVD-KF} implementation is shown to be algebraically equivalent to the conventional KF and, hence, to the earlier published \texttt{SVD-SRKF} variant~\cite{WangSVD1992}.

\begin{rem} \label{rem:2}
Similarly to the \texttt{SVD-SRKF} implementation developed in~\cite{WangSVD1992}, the Cholesky decomposition might be used for the process and measurement noise covariances $\Theta_k$ and $R_k$ in the new \texttt{SVD-KF} designed in this paper. To do that, one needs to replace the matrix product $D^{1/2}_{\Theta_{k-1}}  Q^T_{\Theta_{k-1}}$ in equation~\eqref{svd:5} by an upper-triangular Cholesky factor $\Theta_{k-1}^{1/2}$ (i.e. $\Theta_{k-1} = \Theta_{k-1}^{T/2}\Theta_{k-1}^{1/2}$) and, also, the matrix product  $D^{1/2}_{R_k}  Q^T_{R_k}$ in equations~\eqref{svd:1}, \eqref{svd:2} by an upper-triangular Cholesky factor $R_k^{1/2}$ (i.e. $R_k = R_k^{T/2}R_k^{1/2}$).
\end{rem}

In summary, the \texttt{SVD-KF} implementation provides robust computations meaning that the method works for any covariance matrices $\Theta_k \ge 0$ and $R_k \ge 0$. Meanwhile, the \texttt{SVD-SRKF} counterpart might be unexpectedly interrupted by the error appeared in the function, performing the Cholesky decomposition. From this point of view, the \texttt{SVD-KF} variant might be preferable for solving practical problems. However, the \texttt{SVD-KF} is expected to be slower than the \texttt{SVD-SRKF} because of utilizing that SVD factorization instead of the Cholesky decomposition. If the matrices $\Theta$ and $R$ are time-invariant, then the corresponding factorization is performed only once. This can be done at the initial step, i.e. the corresponding SVD factors of $\Theta$ and $R$ might be pre-computed before the KF recursion. Hence, for constant (over time) covariances $\Theta$ and $R$, the difference in time consumption between the new SVD-KF variants with the SVD and the Cholesky decomposed covariances will be negligible.

Next, the new \texttt{SVD-KF} avoids inversion of $D_{P_{k|k-1}}^{1/2}$ and $D_{P_{k|k}}^{1/2}$ in each iteration step of the filter. As a result, it is expected to be inherently more stable (with respect to roundoff errors) than the previously developed  \texttt{SVD-SRKF} variant~\cite{WangSVD1992}, especially, when the error covariance matrix is close to singular. However, this is achieved by an additional SVD factorization compared to the \texttt{SVD-SRKF}, which makes the new \texttt{SVD-KF} slightly slower. In summary, the novel \texttt{SVD-KF} implementation is expected to be numerically more stable than the \texttt{SVD-SRKF}, but also slower because of that additional SVD factorization.

Furthermore, the new \texttt{SVD-KF} version does not require the triangular matrix inversion, $R_k^{1/2} \in \mathbb{R}^{m\times m}$, compared to the original \texttt{SVD-SRKF} variant. Instead, it involves the $D_{R_{e,k}}^{-1}$ calculation for finding {\it a posteriori} state by equation~\eqref{svd:3}. The mentioned matrix contains the singular values of the innovation covariance $R_{e,k}$.
In future, if possible, we will intend for mitigating this requirement, i.e. to design the SVD-based filter without the mentioned inversion.

Finally, in contrast to the \texttt{SVD-SRKF} published earlier, the new \texttt{SVD-KF} method is conveniently suited for the log likelihood function (log LF) evaluation that is given as follows~\cite{Taylor1986}:
\begin{equation}
{\mathcal L}=-c_0 -
\frac{1}{2} \sum \limits_{k=1}^K \left\{
 \ln\left(\det R_{e,k}\right)+ e_k^T R_{e,k}^{-1}e_k \right\} \label{llf:conv}
\end{equation}
where $c_0 = - (1/2)Km\ln(2\pi)$ is some constant value.

As can be seen, the quantity $R_{e,k}$ required in the log LF computation~\eqref{llf:conv} is not directly available from the original \texttt{SVD-SRKF} implementation via recursion~\eqref{kf:p:X}, \eqref{kf:f:X}, \eqref{svd:p:fac}~-- \eqref{svd:f:K}. Meanwhile, in the new \texttt{SVD-KF} method, the required terms $D_{R_{e,k}}^{1/2}$ and $Q_{R_{e,k}}$ are simply read off from the post-array factors in~\eqref{svd:1}. Formula~\eqref{llf:conv} can be rewritten in terms of the quantities that appear naturally in the \texttt{SVD-KF} as follows:
\begin{equation}
{\mathcal L}\left(Z_1^K\right) =
c_0 - \frac{1}{2} \sum \limits_{k=1}^K \left\{
  \ln\left(\det D_{R_{e,k}}\right)+ \bar e_k^T D^{-1}_{R_{e,k}}\bar e_k
\right\} \label{llf:svd}
\end{equation}
where $\bar e_k$ are the normalized innovations defined as $\bar e_k  =  Q_{R_{e,k}}^Te_k$.

To derive equation~\eqref{llf:svd}, one takes into account the properties of orthogonal matrices and concludes that
  \begin{align}
  \det\left(R_{e,k}\right) & = \det\left(D_{R_{e,k}}\right), &   e_k^T R_{e,k}^{-1} e_k & = \bar e_k^T D^{-1}_{R_{e,k}}\bar e_k. \label{new:formulas}
  \end{align}
By substituting~\eqref{new:formulas} into~\eqref{llf:conv}, we obtain equation~\eqref{llf:svd}.

In conclusion, the novel \texttt{SVD-KF} is expected to be slightly slower than the previously developed \texttt{SVD-SRKF} variant, because of an additional SVD factorization in each iteration step of the \texttt{SVD-KF}. Indeed, there are two SVD in the \texttt{SVD-SRKF} implementation via equations~\eqref{kf:p:X}, \eqref{kf:f:X}, \eqref{svd:p:fac}~-- \eqref{svd:f:K}  against three SVD in the \texttt{SVD-KF} algorithm by formulas~\eqref{kf:p:X},~\eqref{svd:5}~-- \eqref{svd:3} in each recursion step. However, this extra SVD factorization in the \texttt{SVD-KF} variant introduces an important feature into the new implementation that is a simple and convenient way for calculating the log LF.

\section{Numerical Experiments} \label{sec:exp}

In our comparative study we test the new \texttt{SVD-KF} algorithm against  the previously published \texttt{SVD-SRKF} variant~\cite{WangSVD1992} and the conventional KF implementation~\eqref{kf:p:X}~-- \eqref{kf:f:P}. Additionally, we inspect numerically stable SR-based filtering from~\cite[p.~434-436]{KailathSayed2000} (abbreviated further to \texttt{SR-KF}) and the UD-based counterpart~\cite[p.~261]{GrewalAndrews2015} (abbreviated further to \texttt{UD-KF}). The last two techniques are based on the Cholesky decomposition and its `square-root-free' modification of the state error covariance matrix, respectively.

\begin{exmp} \label{ex:1}
 The dynamic  of the in-track motion  of a satellite traveling in a circular orbit is given as follows~\cite[p.~1448]{Rauch1965}:
\begin{align*}
x_{k} & =
\begin{bmatrix}
1 & 1 & 0.5 &  0.5 \\
0 & 1 & 1 & 1 \\
0 & 0 & 1 & 0 \\
0 & 0 & 0 & 0.606
\end{bmatrix}
x_{k-1}\! +\! w_{k-1},
\Theta
=
\begin{bmatrix}
0 & 0 & 0 & 0 \\
0 & 0 & 0 & 0 \\
0 & 0 & 0 & 0 \\
0 & 0 & 0 & q
\end{bmatrix},
 \\
z_k & =
\begin{bmatrix}
1 & 0 & 0 & 0
\end{bmatrix}
x_k + v_k, \quad R = 1, \quad q=0.63 \cdot 10^{-2}
\end{align*}
with zero-mean initial state and $\Pi_0 = diag\{[1, 1, 1, 10^{-2}]\}$.
\end{exmp}

We perform the following set of numerical experiments. Starting from the initial values, the system above is simulated for $k = 1, \ldots, K$, $K=100$, discrete time points in order to generate the ``true'' trajectory of the dynamic state, $x^{exact}_{k}$, $k=1, \ldots, K$, and the corresponding measurements $z_k$, $k=1, \ldots, K$. Next, the inverse problem is solved by the examined filtering methods. More precisely, given the measurement history, each filter under assessment provide the estimated path $\hat x_{k|k}$, $k=1, \ldots, K$ of the unknown dynamic state. For a fair comparative study, the same filter initial values, the same system matrices and the same measurements are passed to all KF variants under examination.
All methods were implemented in the same precision (64-bit
floating point) in MATLAB running on a conventional PC with processor \verb"Intel(R) Core(TM) i5-2410M CPU 2.30 GHz" and with $4$ GB of installed memory (RAM).

The experiment outlined above is repeated for $M=500$ times, i.e. we perform $500$
Monte-Carlo simulations. In each Monte-Carlo run we get ``true'' trajectory of the dynamic state, $x^{exact}_{k}$, $k=1, \ldots, K$, the corresponding measurement history $z_k$, $k=1, \ldots, K$ and the estimated dynamic state path $\hat x_{k|k}$, $k=1, \ldots, K$. Thus, we report the root {\it mean} square error (RMSE) in each component of the state vector, given as follows:
\begin{equation} \label{eq:RMSEx}
\mbox{\rm RMSE}_{x_i}=\sqrt{\frac{1}{MK}\sum \limits_{j=1}^{M} \sum \limits_{k=1}^K \left(x_{i, k}^{j, exact} - \hat x_{i, k|k}^j\right)^2}
\end{equation}
where $M=500$ is the number of Monte-Carlo trials, $K=100$ is the discrete time of the dynamic system, the $x_{i, k}^{j, exact}$ and $\hat x_{i, k|k}^j$ are the $i$-th entry of the ``true'' state vector (simulated) and its estimated value obtained in the $j$-th Monte Carlo trial, respectively. Together with the absolute errors, we calculate the corresponding mean relative errors at each component of the state vector as follows:
\begin{equation} \label{eq:MMRE}
\mbox{\rm MRE}_{x_i} = \sum \limits_{j=1}^{M} \frac{\left|x_{i, K}^{j, exact} - \hat x_{i, K|K}^j\right|}{\left|x_{i, K}^{j, exact}\right|}, \quad i = 1, \ldots, 4.
\end{equation}

Formula~\eqref{eq:MMRE} means that the relative errors are computed only for the steady-state solution, i.e. at the last time point $K=100$. In other words, we exclude the filter transition period (when the filter converges to the steady state) from the relative error computations. Besides, the relative error is possible to compute only for nonzero exact value. Hence, if any entry of the ``true'' state vector tends to zero, then the corresponding relative error does not exist. In our numerical experiments in Example~1, the third component of the ``exact'' state vector (the time path) is always about zero and tends to zero. Instead of the computed relative error, we use '--' for this component. Table~1 summarizes the resulted absolute and relative errors together with the  CPU time (s) averaged over $M=500$ Monte Carlo simulations for each estimator under examination.

\begin{table*}[!t]
\begin{center}
\processtable{The absolute and relative errors with the average CPU time(s) in Example~1, 500 Monte Carlo runs.}
{\tabcolsep6pt\begin{tabular}{rl|cccc|cccc|c}
\hline
Approach & Method & \multicolumn{4}{c|}{RMSE$_{x_i}$} &  \multicolumn{4}{c|}{MRE$_{x_i}$, \%} &  \\
 \hhline{*{2}{|~}*{8}{|-}*{1}{|~}}
  & & $x_1$ & $x_2$ & $x_3$ & $x_4$ & $x_1$ & $x_2$ & $x_3$ & $x_4$  & CPU  \\
\hline
Conventional:             & {\tt KF}       & 0.5762  &  0.3034  & 0.0596  & 0.1073 & 0.0880 & 0.7994 & $-$ & 6.5686  & 0.0078 \\
Cholesky factorization-based:   & {\tt SR-KF}    & 0.5762  &  0.3034  & 0.0596  & 0.1073 & 0.0880 & 0.7994 & $-$ & 6.5686  & 0.0123 \\
UD factorization-based:   & {\tt UD-KF}    & 0.5762  &  0.3034  & 0.0596  & 0.1073 & 0.0880 & 0.7994 & $-$ & 6.5686  & 0.0121 \\
\hline
SVD factorization-based: & {\tt SVD-SRKF} & 0.5762  &  0.3034  & 0.0596  & 0.1073 & 0.0880 & 0.7994 & $-$ & 6.5686  & 0.0139 \\
SVD factorization-based: &{\tt SVD-KF}    & 0.5762  &  0.3034  & 0.0596  & 0.1073 & 0.0880 & 0.7994 & $-$ & 6.5686  & 0.0189 \\
\hline
\end{tabular}}{}
\label{table:1}
\end{center}
\end{table*}

\begin{table*}
\begin{center}
\processtable{The norm of the vector of absolute errors at each component of the dynamic state, i.e. $\|\rm{RMSE}_{x_i}\|_2$, in Example~2, 500 Monte Carlo runs.}
{\tabcolsep6pt\begin{tabular}{l|cccccccccccccc}
\hline
Method & \multicolumn{14}{c}{The $\|\rm{RMSE}_{x_i}\|_2$, $i=1, \ldots, 4$, while growing ill-conditioning $\delta \to \epsilon_{roundoff}$}  \\
\hhline{*{1}{|~}*{14}{|-}}
   & $10^{-1}$ & $10^{-2}$ & $10^{-3}$ & $10^{-4}$ & $10^{-5}$ & $10^{-6}$ & $10^{-7}$ & $10^{-8}$ & $10^{-9}$ &  $10^{-10}$ & $10^{-11}$ & $10^{-12}$ & $10^{-13}$ & $10^{-14}$ \\
\hline
\hline
{\tt KF}       & 0.1230  & 0.0729   & 0.0772  & 0.1002 & 0.0758 & 0.0330 & 0.0519 & \texttt{NaN} & \texttt{NaN} & \texttt{NaN}   & \texttt{NaN}  & \texttt{NaN} & \texttt{NaN}  & \texttt{NaN}\\
{\tt SR-KF}    & 0.1230  & 0.0729   & 0.0772  & 0.1002 & 0.0758 & 0.0330 & 0.0519 & 0.0875       & 0.0593       & 0.0988         & 0.0679        & 0.0588  &  0.0724  & 0.0604\\
{\tt UD-KF}    & 0.1230  & 0.0729   & 0.0772  & 0.1002 & 0.0758 & 0.0330 & 0.0519 & 0.0875       & 0.0593       & 0.0988         & 0.0679        & 0.0588  &  0.0723  & 0.0604\\
\hline
{\tt SVD-SRKF}  &  0.1230  &  0.0729   & 0.0772  & 0.1002 & 0.0758 & 0.0330 & 0.0521 & 0.0905 & \texttt{Inf}   & \texttt{Inf}   & \texttt{NaN} & \texttt{NaN} & \texttt{NaN}    & \texttt{NaN} \\
{\tt SVD-KF}    &  0.1230  &  0.0729   & 0.0772  & 0.1002 & 0.0758 & 0.0330 & 0.0519 & 0.0875 & 0.0593         & 0.0988         & 0.0679       & 0.0588       &  0.0723  & 0.0590\\
\hline
\end{tabular}}{}
\label{table:2}
\end{center}
\end{table*}

Having analyzed the obtained numerical results collected in Table~1, we make a few important conclusions. First, all filtering algorithms explored in this paper provide the same estimation accuracy. This substantiates the theoretical derivation of the new {\tt SVD-KF} implementation presented in Section~\ref{sec:main} and, next, confirms the algebraic equivalence between all KF variants analyzed in this section. We stress that the problem in Example~1 is well posed and well conditioned. As a result, all  KF algorithms produce the same estimates of the dynamic state and work with the same accuracy. To explore the numerical insights of each implementation, Example~2 below represents a set of ill-conditioned problems for the same satellite traveling dynamic. Second, the obtained absolute and relative errors are small for each component of the dynamic state to be estimated, i.e. all KF algorithms assess the unknown state vector accurately.
Third, we observe that the conventional KF implementation works faster than any other its counterpart. The SVD-based filters are slower than the SR-based or UD-based algorithms, because the SVD factorization is more  computationally expensive than the Cholesky decomposition. Finally, among the examined SVD-based methods, the previously published {\tt SVD-SRKF} variant~\cite{WangSVD1992} is the fastest. This conclusion is reasonable because the {\tt SVD-SRKF} requires only two SVD factorizations in each iteration step while the novel {\tt SVD-KF} method implies three SVD factorizations. The increased execution time is the price that we pay for the log LF computation feature appended into our newly-designed SVD-based filter. The difference in CPU time consumed by the original {\tt SVD-SRKF} algorithm and the new {\tt SVD-KF} method is about $0.0042$ (s) in Example~1, that is the {\tt SVD-SRKF} is about $35$\% faster in comparison to the {\tt SVD-KF}.

Next, we wish to explore the numerical insights of the SVD-based filters discussed in this paper. For that, we equip Example~\ref{ex:1} by ill-conditioned measurement scheme as explained in~\cite{GrewalAndrews2015}.
\begin{exmp} \label{ex:2}
 Consider a satellite dynamic in Example~\ref{ex:1} with
\begin{align*}
z_k & =
\begin{bmatrix}
1 & 1 & 1 & 1\\
1 & 1 & 1 & 1+\delta
\end{bmatrix}
x_k + v_k, \quad v_k \sim {\cal N}(0, R)
\end{align*}
where  $R = \delta^2 I_2$  and  the initial conditions $x_0 \sim {\cal N}(0, \Pi_0)$, $\Pi_0 = I_4$. To simulate roundoff we assume that
$\delta^2<\epsilon_{roundoff}$, but $\delta>\epsilon_{roundoff}$
where $\epsilon_{roundoff}$ denotes the unit roundoff
error\footnote{Computer roundoff for floating-point arithmetic is
often characterized by a single parameter $\epsilon_{roundoff}$,
defined as the largest number such that either
$1+\epsilon_{roundoff} = 1$ or $1+\epsilon_{roundoff}/2 = 1$ in
machine precision. }.
\end{exmp}

As in the previous set of numerical experiments, all KF implementations were implemented in the same precision (64-bit floating
point) in MATLAB where the unit roundoff error is $2^{-53}\approx
1.11 \cdot 10^{-16}$. The MATLAB function \verb"eps" is twice the unit
roundoff error. Hence, we repeat the set of numerical experiments outlined above for various values of the ill-conditioning parameter $\delta$ such that $\delta \to \epsilon_{roundoff}$, i.e. $\delta$ tends to the machine precision limit.
 This allows for observing the discrete-time algebraic Riccati equation degradation because the matrix inversion in~\eqref{kf:f:K} will be applied to an ill-conditioned matrix. More precisely, this situation corresponds to the third reason of ill-conditioning listed in \cite[p.~288]{GrewalAndrews2015}. The obtained numerical results are summarized in Table~2.

Having analyzed the outcomes collected in Table~2, we conclude that the conventional KF technique possesses the worst performance among the examined filter implementations. When $\delta$ is large, which corresponds to the well-conditioned problems, all estimators work well and with the same estimation accuracy. However, while growing ill-conditioning, the conventional KF degrades faster than any other algorithm under examination. For $\delta = 10^{-8}$ it produces \verb"NaN" that means 'Not a Number' in MATLAB and this is the sign of a failure of any computational technique.

Next, we compare the factored-form KF implementations and observe  that the previously published \texttt{SVD-SRKF} variant performs markedly worse than any other factorized-form KF under assessment. Indeed, the {\tt SVD-SRKF} method fails when $\delta = 10^{-9}$, meanwhile all other factored-form KF implementations manage this ill-conditioned situation. They work accurately, i.e. with small estimation errors, until very small ill-conditioning parameter values. Thus, our theoretical expectations presented in Sections~3 and 4 are realized. More precisely, the earlier published \texttt{SVD-SRKF} performs faster than the new \texttt{SVD-KF}, but less accurately compared not only to the novel \texttt{SVD-KF}, but also in comparison to all other factored-form KF implementations. Finally, the SVD-based filtering via the new \texttt{SVD-KF} variant is slightly more accurate than the SR- and the UD-based approaches. This conclusion is in line with our recent analysis in~\cite{TsyganovaOnline} where we have extended the advantageous \texttt{SVD-KF}, developed in this paper, on the corresponding filter derivative computations required in {\it parameter estimation} problems.

\section{Concluding remarks}

In this paper, the new SVD factorization-based KF implementation is derived. It is shown to outperform the previously proposed SVD-based KF variant in estimation accuracy when solving ill-conditioned state estimation problem. Our preliminary analysis based on ill-conditioned tests suggest that the novel SVD-based KF approach seems to be slightly more accurate than the SR Cholesky- and UD-based filtering strategies. However, a rigorous theoretical analysis of the error propagation (due to numerical roundoff) should be performed to explore the numerical insights of the mentioned algorithms. This is an open question for a
future research. Another open problem is a development of {\it information}-type SVD-based KF implementations that process the inverse of the error covariance.

Finally, the results of our numerical experiments suggest that the previously proposed \texttt{SVD-SRKF} algorithm performs markedly worse than any other factorized-form KF implementation under assessment, i.e. it degrades faster than the SR-based, the UD-based and the new SVD-based methods in line with the growing problem's ill-conditioning. There are two main sources of its numerical instability: i) the use of equation~\eqref{kf:P:eq1} for computing {\it a posteriori} error covariance matrix instead of the symmetric Joseph stabilized form~\eqref{kf:P:eq3}, and ii) the requirement for inversion of the diagonal matrix containing the square roots of  singular values of the error covariance matrix in each iteration step. For a moment, it is difficult to distinguish and to assess the contribution and impact of each particular factor to the numerical instability of the \texttt{SVD-SRKF}, separately. This is an interesting issue for a future research.


\section*{Acknowledgments}

The first author thanks the support of Portuguese National Fund ({\it Funda\c{c}\~{a}o para a
Ci\^{e}ncia e a Tecnologia}) within the scope of project UID/Multi/04621/2013.

\bibliographystyle{iet}


\end{document}